\documentclass[reqno, a4paper, 10pt]{amsart}  

\usepackage{epsfig} 
\usepackage{amsmath} 
\usepackage{amssymb}  
\usepackage{accents}
\usepackage{pstricks}
\usepackage{pstricks-add}
\usepackage{url}
\usepackage{bbm}

\title[Converse Lyapunov Theorems in Banach Spaces]{Converse Lyapunov Theorems for Switched Systems in Banach and Hilbert Spaces}

\author{Falk M. Hante and Mario Sigalotti}

\thanks{This work was supported by
the ANR project ArHyCo, Programme 
ARPEGE, contract number ANR-2008 SEGI 004  01-30011459.\\
F. M. Hante and M. Sigalotti are with Institut
\'Elie Cartan, UMR 7502, 
 BP 239, Vand\oe uvre-l\`es-Nancy 54506,  France
 and CORIDA, INRIA Nancy -- Grand Est, France. 
E-Mail: \url{hante@iecn.u-nancy.fr} and \url{mario.sigalotti@inria.fr}.
}

\newcommand{\Lcal}{\mathcal{L}}
\newcommand{\B}{\mathcal{B}}
\newcommand{\Scal}{\mathcal{S}}
\renewcommand{\:}{\mathcal{\colon}}
\newcommand{\NN}{\mathbb{N}}
\newcommand{\QQ}{\mathbb{Q}}
\newcommand{\RR}{\mathbb{R}}
\newcommand{\One}{\mathbbm{1}}
\newcommand{\Lsup}{\bar{L}}
\newcommand{\Linf}{\underaccent{\bar}{L}}

\newcommand{\wot}{\stackrel{\mbox{{\tiny{WOT}}}}{\longrightarrow}}
\newcommand{\la}{\langle}
\newcommand{\ra}{\rangle}

\DeclareMathOperator{\argmax}{arg\,max}

\theoremstyle{plain}
\newtheorem{theorem}{Theorem}

\newtheorem{corollary}[theorem]{Corollary}
\newtheorem{lemma}[theorem]{Lemma}
\theoremstyle{remark}
\newtheorem{remark}[theorem]{Remark}
\newtheorem{example}[theorem]{Example}

\begin{document}

\maketitle

\begin{abstract}
We consider switched systems on 
Banach and Hilbert spaces  
governed by
strongly continuous one-parameter semigroups of linear evolution operators. 
We 
provide necessary and sufficient conditions 
for their global exponential stability, uniform with respect to  
the switching signal, 
in terms of the existence of a 
Lyapunov function common to all modes.
\end{abstract}

\section{INTRODUCTION}
It is well known that the existence of a \emph{common Lyapunov function}
is necessary and sufficient for the global uniform
asymptotic
stability of 
finite-dimensional continuous-time switched 
dynamical systems 
\cite{LinSontagWang1996}. 
In the \emph{linear} finite-dimensional case, 
the existence of a common Lyapunov function
is actually equivalent to 
global uniform \emph{exponential} 
stability \cite{MolchanovPyatnitskiy1986,MolchanovPyatnitskiy1989} and, provided that the system has finitely many  modes, 
the Lyapunov function can be taken polyhedral or polynomial (see 
\cite{Blanchini1994,BlanchiniMiani1999,DayawansaMartin1999} and also 
\cite{BraytonTong1980}
for a discrete-time version).  
A special role in the switched control literature has been played by common
quadratic Lyapunov functions, since their existence can be tested rather efficiently (see the surveys \cite{LinAntsaklis2009,ShortenEtAl2007} and the references therein). It is known, however, that 
the 
existence of a common quadratic Lyapunov function
is not necessary for the 
global uniform exponential 
stability of a linear switched system with finitely many modes. Moreover, there exists no uniform upper bound on the minimal  
degree of a common polynomial Lyapunov function \cite{MasonBoscainChitour2006}. 

The scope of this paper is to 
prove that 
the existence of a {common Lyapunov function}
is equivalent to the global uniform exponential 
stability of 
in\-fi\-ni\-te-di\-men\-sio\-nal switched systems of the type
\begin{equation}\label{eq:ssOpForm}
\left\{\begin{aligned}
 \frac{d}{dt}\,x(t)	&= A_{\sigma(t)} x(t),\quad t > 0\\
 	x(0)		&= x \in X
\end{aligned}\right.
\end{equation}
where each $A_j$ is a (possibly unbounded) operator generating a strongly continuous
semigroup $T_j(t)$ on a Banach space $X$ and $\sigma(\cdot)$ belongs to the class of 
piecewise constant switching signals with values in an index set $Q$.

Such systems provide a convenient design paradigm for modeling a wide variety of complex
processes comprising distributed parameters, see \cite{HanteLeugeringSeidman2009} and the references therein
for examples in the context of networked transport systems.

Except for special cases, 
that is,
when $X$ has a Hilbert structure and the infinitesimal generators 
$A_j$ commute pairwise \cite{Sasane2005}, when the switching signals satisfy a 
dwell-time constraint \cite{MichelSunMolchanov2005} or when $A_j$ is a linear convection-reaction operator with reflecting 
boundary conditions \cite{AminHanteBayen2008}, global uniform exponential stability of systems such 
as \eqref{eq:ssOpForm} has not been investigated (up to our knowledge).

The characterization of exponential stability for a single linear dynamical system on 
Banach and Hilbert spaces dates back to Datko \cite{Datko1968} and Pazy \cite{Pazy1972} and has, 
since then, seen a broad range of applications in control theory for partial differential equations (see, for instance, \cite{TucsnakWeiss2009}).
However, we recall that exponential stability of all subsystems (with $\sigma(t) \equiv j$ 
fixed in \eqref{eq:ssOpForm}) is of course necessary but not sufficient for the global uniform 
exponential stability with respect to all possible 
switching laws $\sigma(\cdot)$. This is a classical result for the finite dimensional case and we give an
infinite dimensional variant with interesting destabilizing properties
in Example~\ref{ex:notGUES} below.

Our starting point will be a switching system of the general form
\begin{equation}\label{eq:ss:intro}
\left\{\begin{aligned}
 x(t_{k+1}) &= T_{\sigma(t_k)}(t_{k+1}-t_k) x(t_k),\quad k\in\mathbb{N},\\
 x(0)&=x_0 \in X, 
 \end{aligned}\right.
\end{equation}
where $\sigma\: [0,\infty) \to Q$ is a piecewise constant right-continuous  
switching signal with switching times $0=t_0<t_1<\dots<t_k<\cdots$. 
Each $t\mapsto T_j(t)$, $j\in Q$, is a strongly continuous semigroup on a Banach space $X$.
If $t\in(t_k,t_{k+1})$, then $x(t)=T_{\sigma(t_k)}(t-t_k) x(t_k)$. 
This framework includes, in particular, 
switched dynamical systems such as \eqref{eq:ssOpForm}, even 
with the infinitesimal generators $A_j$ not sharing a common domain and also in the case of 
infinitely many available modes $Q$. (The 
semigroup formulation \eqref{eq:ss:intro} 
corresponds to the choice of mild 
solutions for the Cauchy problem \eqref{eq:ssOpForm}. The two formulations are clearly equivalent when $X$ has finite dimension, with $T_j(t)=e^{t A_j}$.)

The main result of this paper is 
that
the following three conditions are equivalent:
\emph{
\begin{itemize}
 \item[(A)]
There exist two constants $K \geq 1$ and $\mu>0$
such that, 
for every $\sigma(\cdot)$ and every $x_0$,  the solution $x(\cdot)$ to \eqref{eq:ss:intro} satisfies
\begin{equation}\label{thm0:GUES}
 \|x(t)\|_{X} \leq Ke^{-\mu t}\|x_0\|_{X}, \quad t \geq 0.
\end{equation}
\item[(B)] There exist two constants $M \geq 1$ and $\omega>0$
such that, 
for every $\sigma(\cdot)$ and every $x_0$,  the solution $x(\cdot)$ to \eqref{eq:ss:intro} satisfies
\begin{equation}\label{thm0:UEBound}
 \|x(t)\|_{X} \leq Me^{\omega t}\|x_0\|_{X}, \quad t \geq 0,
\end{equation}
and there exists $V\: X \to [0,\infty)$
 such that 
 $\sqrt{V(\cdot)}$ is a norm on $X$, 
	\begin{equation}\label{thm0:VBound}
 		V(x) \leq C \|x\|_X^2, \quad x\in X
	\end{equation}
	for a constant $C>0$ and
	\begin{equation}\label{thm0:LjVest}
 		\liminf_{t \downarrow 0} \frac{V(T_j(t)x) - V(x)}{t} \leq -\|x\|_X^2,\quad
		j \in Q,~x\in X.
	\end{equation}
\item[(C)] There exists $V\: X \to [0,\infty)$
 such that 
 $\sqrt{V(\cdot)}$ is a norm on $X$, 
	\begin{equation}\label{thm00:VBound}
 		c\|x\|_X^2\leq V(x) \leq C \|x\|_X^2, \quad x\in X
	\end{equation}
	for some constants $c,C>0$ and
	\begin{equation}\label{thm00:LjVest}
 		\liminf_{t \downarrow 0} \frac{V(T_j(t)x) - V(x)}{t} \leq -\|x\|_X^2,\quad j \in Q,~x\in X.
	\end{equation}
		\end{itemize}
}

The equivalence between (A) and (C) extends to infinite-dimensional systems the well-known result obtained in \cite{MolchanovPyatnitskiy1986}
in the finite-dimensional setting.

Conditions \eqref{thm0:VBound} and \eqref{thm00:VBound} are redundant in the case of finite-dimensional systems, since 
$\sqrt{V(\cdot)}$ and $\|\cdot\|_X$ are comparable, by compactness of the unit sphere. 
Hence, condition \eqref{thm0:UEBound} in (B) could be dropped for finite-dimensional systems. This is not the case 
for infinite-dimensional ones, as illustrated in 
Remark~\ref{rem:Thm1sharp} by an example.

From the point of view of applications, 
condition (B), 
imposing less conditions on $V$ than (C), is better suited for 
establishing that (A) holds (although the uniform exponential growth boundedness needs also be proved).
On the other hand, the implication (A)$\Rightarrow$(C)  can be used to 
select a Lyapunov function with tighter requirements.

The construction of a common Lyapunov function
satisfying (B), under the assumption that (A) holds true, follows the same lines as in finite dimension. In particular, a possible choice of the Lyapunov function is 
\begin{equation*}
V(x_0)=\sup\left\{ \int_0^\infty \|x(t)\|^2 dt : x(\cdot)~\text{solution to \eqref{eq:ss:intro} for some}~\sigma\right\}.
\end{equation*}
(Alternatively, one could take $V(x_0)=\int_0^\infty \sup_{\sigma(\cdot)}\|x(t)\|^2 dt$, as done in 
\cite{HanteSigalotti2010}.)

The construction of a Lyapunov function
satisfying (C) is similar, but one has to 
augment \eqref{eq:ss:intro} with a further mode $T_{j^*}(t)=e^{-\mu t}I$ (where $\mu>0$ is the constant appearing in (A) and $I$ denotes the identity on $X$) 
and to consider all the solutions to this augmented system in the definition of $V$.

In the case of an exponentially stable single mode ($Q=\{0\}$), it was observed by Pazy \cite{Pazy1972}
that $x\mapsto \int_0^\infty \|T_0(t)x\|^2\,dt$ defines a Lyapunov function that is comparable with the squared norm if and only if $T_0$ extends to an exponentially stable strongly continuous group. Notice that, as a consequence of the implication (A)$\Rightarrow$(C), even if $T_0$ does not admit an extension to a group, a Lyapunov function comparable with the squared norm can still be found (see Remark~\ref{one-mode-L}).




Concerning the regularity of the Lyapunov functions obtained through the construction described above, they are always convex and continuous (since $\sqrt{V(\cdot)}$ is a norm).
 In the special case in which $X$ is a separable Hilbert space, we 
also prove the Fr\'echet directional differentiability of $V$ and we establish a characterization of the directional Fr\'echet derivatives. 

The paper is organized as follows. 
In Section~\ref{sec:notations} we introduce the main notations and discuss a motivation example. 
Section~\ref{sec:first}
provides a first necessary and sufficient condition for global uniform exponential stability
in terms of the existence of a common Lyapunov function, namely, the equivalence of  (A) and (B)
(Theorem~\ref{thm:VnonComparable}). 
We also discuss the possible redundancies of 
condition (B),
showing that 
\eqref{thm0:UEBound}
cannot be removed (Remark~\ref{rem:Thm1sharp}). 
In Section~\ref{sec:second} a second converse Lyapunov theorem is proved, establishing that (A) and (C) are equivalent (Theorem~\ref{thm:VComparable}).
In Section~\ref{sec:Hilbert} we show the Fr\'echet differentiability of the common Lyapunov functions constructed in the previous sections when $X$ is a separable Hilbert space (Corollary~\ref{corl:Hilbert}).
In Section~\ref{sec:conclusions} we give some final remarks and point to open problems.

\section{
NOTATIONS AND PRELIMINARIES}\label{sec:notations} 

By $\NN$, $\QQ$ and $\RR$ we denote the set of natural, rational and real numbers, respectively. Further, let $X$ be a Banach space, $\Lcal(X)$ be the space of bounded linear operators on $X$, $Q$ be a countable set 
and, for all $j \in Q$, let $t\mapsto T_j(t)\in \Lcal(X)$, $t \geq 0$ be a strongly continuous semigroup.

We wish to investigate the qualitative behavior of 
\begin{equation}\label{eq:ss}
 x(t) = T_{\sigma(\cdot)}(t)x 
\end{equation}
for $x\in X$, where $\sigma\: [0,\infty) \to Q$ is a piecewise constant 
switching signal and
\begin{equation}\label{eq:TsigmaDef}
T_{\sigma(\cdot)}(t)=T_{j_p}(t-\tau_p)T_{j_{p-1}}(\tau_p-\tau_{p-1}) \cdots T_{j_{1}}(\tau_{1})
\end{equation}
for $\sigma(\cdot)$ equal to $j_k$ on $(\tau_{k-1},\tau_k)$ for $k=1,\dots,p+1$ and 
\begin{equation*}
0=\tau_0<\tau_1<\dots<\tau_{p+1}=t.
\end{equation*}
In particular, we wish to study the asymptotic behavior of $x(t)$ as $t$ tends 
to $+\infty$, uniformly with respect to the switching law $\sigma(\cdot)$ in the set $\Sigma$ of all piecewise constant switching signals. 
We note that, for any given $\sigma(\cdot)\in\Sigma$, the operator
$T_{\sigma(\cdot)}(t)\in\Lcal(X)$ is strongly continuous with respect to $t$ 
($\lim_{t\downarrow t_0}\|T_{\sigma(\cdot)}(t)-T_{\sigma(\cdot)}(t_0)\|=0$)
and satisfies 
\begin{equation}\label{eq:SemigoupProperty}
 T_{\sigma(\cdot)}(t+s)=T_{\sigma_s(\cdot)}(t)T_{\sigma(\cdot)}(s)
\end{equation}
for some switching signal $\sigma_s(\cdot)\in\Sigma$ depending on $s$, but in general,
$T_{\sigma(\cdot)}(t)$ does not satisfy the semigroup property, i.\,e., 
equation \eqref{eq:SemigoupProperty} with $\sigma_s(\cdot)$ replaced by $\sigma(\cdot)$ (independently of $s$).

For a function $V\: X \to [0,\infty)$ we define the generalized derivative
\begin{equation}\label{eq:LinfjVdef}
 \Linf_j V(x)=\liminf_{t \downarrow 0} \frac{V(T_j(t)x) - V(x)}{t},
\end{equation}
noting the possibility that $|\Linf_j V(x)|=\infty$ for some $x\in X$ and $j \in Q$.
Further, we call a switched system \eqref{eq:ss} (completely determined
by $\{T_{j}\}_{j \in Q}$) \emph{globally uniformly exponentially stable} when there exist
constants $K \geq 1$ and $\mu>0$ such that
\begin{equation}\label{def:GUES}
 \|T_{\sigma(\cdot)}(t)\|_{\Lcal(X)} \leq K e^{-\mu t}, \quad t \geq 0,~\sigma(\cdot)\text{-uniformly}.
\end{equation}
It is clear that \eqref{def:GUES} implies
\begin{equation*}
\|x(t)\|_X \leq K e^{-\mu t} \|x\|_X, \quad t \geq 0
\end{equation*}
globally for all $x \in X$ and uniformly for all $\sigma(\cdot)\in\Sigma$
justifying the terminology. We point out that \eqref{def:GUES} implies
\emph{strong attractivity at the origin}, i.\,e.,
\begin{equation}\label{def:strongStab}
 \lim_{t\to\infty}\|T_{\sigma(\cdot)}(t)x\|_X = 0, \quad x \in X,~\sigma(\cdot)\in\Sigma, 
\end{equation}
and \emph{uniform stability}, i.\,e., 
\begin{equation}\label{def:unifStab}
\left\{\begin{aligned}
&\text{for all}~\varepsilon>0~\text{there exists a}~\delta>0,\\
&\text{independent of}~\sigma(\cdot),~\text{such that}~\|x\|_X<\delta~\text{implies}\\
&\quad \|T_{\sigma(\cdot)}(t)x\|_{X}<\varepsilon,~t \geq 0,~\sigma(\cdot)\text{-uniformly}
\end{aligned}\right.
\end{equation}
but that the converse implication is
false in general, even for a single mode: As a counterexample it suffices to take the left translation semigroup 
defined by 
\begin{equation*}
(T(t)f)(s):=f(s+t) 
\end{equation*}
on the Lebesgue space $X=L^1(\RR_+)$. This is in contrast to the equivalence of \eqref{def:strongStab} and \eqref{def:GUES} when $X$ is a $n$-dimensional real coordinate space, $Q$ is finite, and $T_{j}(t)$ is given by the matrix exponential $e^{A_jt}$ for some real $n \times n$-matrix
$A_j$, as a consequence of Fenichel's Uniformity Lemma (see, for instance, \cite[\S 5.2]{ColoniusKliemann2000}).

Before turning our attention to necessary and sufficient conditions for global
uniform exponential stability, we give an example of a switched system exhibiting
illustrative instability properties, though the subsystems are exponentially stable.

\begin{figure*}[t]
\centering
\begin{pspicture}(12,10)
\put(0,4.5){\scalebox{1}[1]{
\begin{pspicture}(7,5)
\pspolygon[linestyle=none,fillstyle=crosshatch*,hatchangle=45,hatchwidth=0.1pt](2.5,1)(1,2.5)(1,4)(2.5,2.5)(2.5,1)
\pspolygon[linestyle=none,fillstyle=vlines,hatchangle=45,hatchwidth=0.1pt](1,1)(2.5,1)(1,2.5)(1,1)
\pspolygon[linestyle=none,fillstyle=vlines,hatchangle=45,hatchwidth=0.1pt](2.5,1)(4,1)(2.5,2.5)(2.5,1)
\psline{->}(1,1)(4.5,1)
\psline[linestyle=dotted]{-}(1,1)(1,4.5)
\psline{->}(2.5,1)(2.5,4.5)
\psline[linestyle=dotted]{-}(4,1)(4,4.5)
\psline[arrowsize=0.2,linewidth=0.5pt]{->}(4,1)(1,4)
\psline[arrowsize=0.2,linewidth=0.5pt]{->}(2.5,1)(1,2.5)
\rput(0.5,2.5){$t$}
\rput(5,0.75){$s$}
\psline[linestyle=dotted]{-}(0.95,4)(4,4)
\psline{-}(0.95,4)(1.05,4)
\rput(0.75,4){$2$}
\psline{-}(1,0.95)(1,1.05)
\rput(0.9,0.75){$-1$}
\psline{-}(4,0.95)(4,1.05)
\rput(4,0.75){$1$}
\end{pspicture}
}}

\put(0,0){\scalebox{1}[1]{
\begin{pspicture}(7,5)
\pspolygon[linestyle=none,fillstyle=crosshatch,hatchangle=45,hatchwidth=0.1pt](2.5,1)(4,2.5)(4,4)(2.5,2.5)(2.5,1)
\pspolygon[linestyle=none,fillstyle=vlines,hatchangle=-45,hatchwidth=0.1pt](2.5,1)(4,1)(4,2.5)(2.5,1)
\pspolygon[linestyle=none,fillstyle=vlines,hatchangle=-45,hatchwidth=0.1pt](2.5,1)(2.5,2.5)(1,1)(2.5,1)
\psline{->}(1,1)(4.5,1)
\psline[linestyle=dotted]{-}(1,1)(1,4.5)
\psline{->}(2.5,1)(2.5,4.5)
\psline[linestyle=dotted]{-}(4,1)(4,4.5)
\psline[arrowsize=0.2,linewidth=0.5pt]{->}(1,1)(4,4)
\psline[arrowsize=0.2,linewidth=0.5pt]{->}(2.5,1)(4,2.5)
\rput(0.5,2.5){$t$}
\rput(5,0.75){$s$}
\psline[linestyle=dotted]{-}(0.95,4)(4,4)
\psline{-}(0.95,4)(1.05,4)
\rput(0.75,4){$2$}
\psline{-}(1,0.95)(1,1.05)
\rput(0.9,0.75){$-1$}
\psline{-}(4,0.95)(4,1.05)
\rput(4,0.75){$1$}
\end{pspicture}
}}

\put(6.5,0){\scalebox{1}[1]{
\begin{pspicture}(6,10)
\psline{->}(1,1)(5,1)
\psline{->}(1,1)(1,9)

\rput(-0.2,4.5){$\|T_{\sigma(\cdot)}(t)\|_{\Lcal(X)}$}
\rput(5.5,0.6){$t$}

\psline{-}(0.95,1.25)(1.05,1.25)
\rput(0.8,1.25){$1$}

\psline{-}(4.5,0.95)(4.5,1.05)
\rput(4.5,0.7){$1$}

\psline[linewidth=0.5pt]{-}(1,1.5)(1.5,1.5)
\psline[linestyle=dotted]{-}(1.5,1)(1.5,2)
\psline[linewidth=0.5pt]{-}(1.5,2)(2,2)
\psline[linestyle=dotted]{-}(2,1)(2,3)
\psline[linewidth=0.5pt]{-}(2,3)(2.5,3)
\psline[linestyle=dotted]{-}(2.5,1)(2.5,5)
\psline[linewidth=0.5pt]{-}(2.5,5)(3,5)
\psline[linestyle=dotted]{-}(3,1)(3,9)
\psline[linewidth=0.5pt]{-}(3,9)(3.5,9)
\psline[linestyle=dotted]{-}(3.5,1)(3.5,9)

\psline{-}(2,0.65)(2,0.75)
\psline{-}(2.5,0.65)(2.5,0.75)
\psline{-}(2,0.7)(2.1,0.7)
\psline{-}(2.4,0.7)(2.5,0.7)
\rput(2.25,0.7){$\delta$}
\end{pspicture}
}}
\end{pspicture}
\caption{Illustration of $T_j(t)$, $j=1,2$, (left) and the blow up of the operator norm of $T_{\sigma(\cdot)}(t)$ (right) in Example~\ref{ex:notGUES}.\label{f::1}}
\end{figure*}
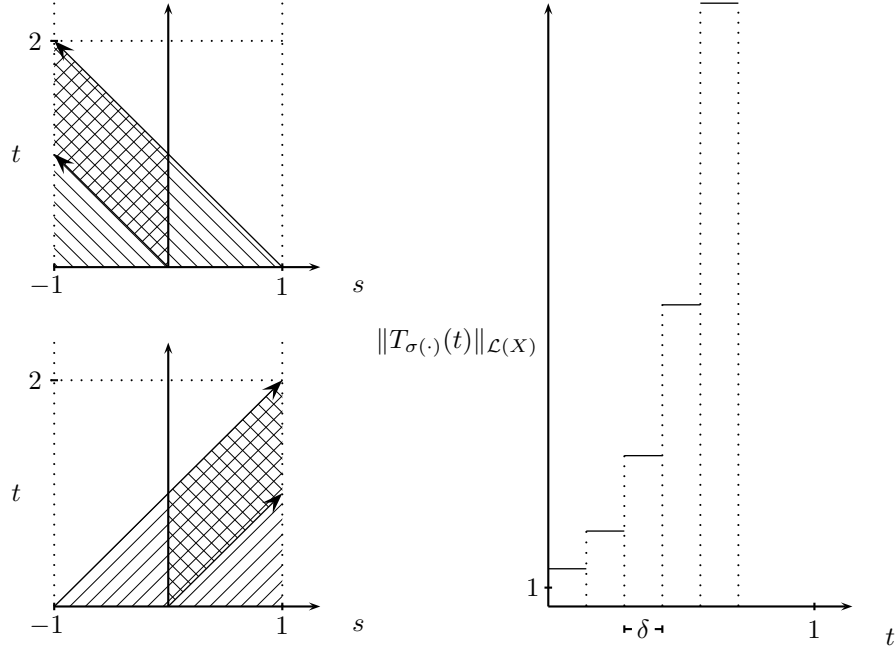

\begin{example}\label{ex:notGUES} Consider the bimodal system $\{T_j(t)\}_{j=1,2}$ with
$T_j(t)$ defined on the Lebesgue space $X=L^1(-1,1)$ by
\begin{equation*}
 \left(T_1(t)f\right)(s)=\begin{cases}
                         	2 f(s+t),	& s \in [-1,1-t]\cap[-t,0]\\
                         	f(s+t),		& s \in [-1,1-t]\setminus[-t,0]\\
				0,		& s \in (1-t,1]
                         \end{cases}
\end{equation*}
and
\begin{equation*}
 \left(T_2(t)f\right)(s)=\begin{cases}
                         	2 f(s-t),	& s \in [-1+t,1]\cap[0,t]\\
                         	f(s-t),		& s \in [-1+t,1]\setminus[0,t]\\
				0,		& s \in [-1,-1+t).
                         \end{cases}
\end{equation*}
Notice that both $T_1(\cdot)$ and $T_2(\cdot)$ are nilpotent semigroups, since 
$T_1(t)=T_2(t)=0$ for $t\geq 2$. In particular, each of them is exponentially stable.

It is easy to see that for
suitable
switching signals $\sigma(\cdot)\in\Sigma$, e.\,g., switching at $\tau_k=k\delta$, $k\in\NN$, 
for a fixed $\delta<1$, 
\begin{equation*}
\|T_{\sigma(\cdot)}(t)\|_{\Lcal(X)}\to+\infty \quad\text{as}~t \to +\infty.
\end{equation*}
In fact, the speed of blow-up is not uniformly exponentially bounded over the set of all possible $\sigma(\cdot)$, i.\,e.,
for any fixed $t>0$, we have 
\begin{equation*}
\|T_{\sigma(\cdot)}(t)\|_{\Lcal(X)}\geq 2^{\lceil\frac{t}{\delta}\rceil}\to+\infty \quad\text{as}~\delta \to 0,
\end{equation*}
with $\lceil\tau\rceil=\min\{k\in \mathbb{N} : \tau\leq k \}$ for $\tau>0$ (see Figure~\ref{f::1}). This can be seen by taking $L^1(-1,1)$-functions $f$ of norm one, identically constant near $x=0$ on progressively smaller intervals and zero elsewhere. \hfill$\square$
\end{example}

\section{FIRST CONVERSE LYAPUNOV THEOREM}\label{sec:first}
In this section we establish a first equivalence result for the global uniform exponential stability of a switched of the form~\eqref{eq:ss}. The crucial step is given by the following lemma, related to the blow-up phenomenon illustrated in Example~\ref{ex:notGUES} in the previous section. It 
is a variant of a result obtained in \cite{Triggiani1994} in the framework of 
strongly continuous semigroups. While extending the property to 
switched system of the form \eqref{eq:ss}, the proof given in \cite{Triggiani1994}
should be modified in order to replace 
the semigroup property by \eqref{eq:SemigoupProperty}.
We include the modified proof for the sake of completeness.
\begin{lemma}\label{lem:IntUBound} Assume that 
\begin{itemize}
\item[i)] there exist constants $M \geq 1$ and $\omega>0$
such that
\begin{equation*}
 \|T_{\sigma(\cdot)}(t)\|_{\Lcal(X)} \leq Me^{\omega t}, \quad t \geq 0,~\sigma(\cdot)\text{-uniformly};
\end{equation*}
\item[ii)] there exists a constant $C>0$ 
and some $p\in[1,\infty)$
such that
\begin{equation*}
 \int_0^\infty \|T_{\sigma(\cdot)}(t)x\|^p_X\,dt \leq C\|x\|^p_X, \quad x \in X,~\sigma(\cdot)\text{-uniformly}.
\end{equation*}
\end{itemize}
Then, there exist constants $K \geq 1$ and $\mu>0$ 
such that
\begin{equation*}
 \|T_{\sigma(\cdot)}(t)\|_{\Lcal(X)} \leq Ke^{-\mu t}, \quad t \geq 0,~\sigma(\cdot)\text{-uniformly}.
\end{equation*}
\mbox{}
\end{lemma}
\begin{proof}
First, we show that under the assumptions i) and ii), for every $x \in X$, there exists a constant $C_x>0$ 
such that
\begin{equation}\label{eq:XNormUBound}
 \|T_{\sigma(\cdot)}(t)x\|_{X} \leq C_x,~t \geq 0, \quad\sigma(\cdot)\text{-uniformly}
\end{equation}
and that, for all $\sigma(\cdot)$ and for all $x \in X$,
\begin{equation}\label{eq:StrongStab}
 \lim_{t\to+\infty}\|T_{\sigma(\cdot)}(t)x\|_{X}=0.
\end{equation}
To this end, let 
$t>\frac{1}{\omega}$ and set $\Delta(t)=[t-\frac{1}{\omega},t]$.
Observe that, for every $\sigma(\cdot)$ and every $\tau\in\Delta(t)$, there exists a
$\sigma_{\tau}(\cdot)$ such that
\begin{equation}\label{eq:DeltaiEst}
\begin{split}
 \|T_{\sigma(\cdot)}(t)x\|_X	&=\|T_{\sigma_{\tau}(\cdot)}(t-\tau)T_{\sigma(\cdot)}(\tau)x\|_X\\
				&\leq \|T_{\sigma_{\tau}(\cdot)}(t-\tau)\|_{\Lcal(X)}\|T_{\sigma(\cdot)}(\tau)x\|_X.
\end{split}
\end{equation}
Moreover, by assumption i)
and by definition of $\Delta(t)$, we have
\begin{equation}\label{eq:OpNormOnDeltaiEst}
 \|T_{\sigma_{\tau}(\cdot)}(t-\tau)\|_{\Lcal(X)} \leq M e^{\omega(t-\tau)}\leq M e^{\omega\frac{1}{\omega}}=Me,
\end{equation}
yielding
\begin{equation}\label{eq:XY}
 \|T_{\sigma(\cdot)}(\tau)x\|_X \geq \frac{\|T_{\sigma(\cdot)}(t)x\|_X}{Me}, \quad\tau\in\Delta(t).
\end{equation}

Now suppose \eqref{eq:XNormUBound} does not hold.
Then, there exist $x \in X$, a sequence of switching signals $(\sigma_i(\cdot))_{i\in\mathbb{N}}$ in $\Sigma$ and a sequence of times $(t_i)_{i\in\mathbb{N}}$
such that 
\begin{equation}\label{eq:XNormDiv}
 \delta_i=\|T_{\sigma_i(\cdot)}(t_i)x\|_X \to +\infty \quad\text{as}~i\to+\infty.
\end{equation}
Assumption i) guarantees that $t_i$ is diverging. Without loss of generality, $t_i>\frac1\omega$ for every $i\in\mathbb{N}$. 

For $\tau\in\Delta(t_i)$, \eqref{eq:XY} yields
\begin{equation}\label{eq:XNormOnDeltaiEst}
 \|T_{\sigma_i(\cdot)}(\tau)x\|_X \geq \frac{\delta_i}{Me}, \quad\tau\in\Delta(t_i).
\end{equation}
Hence, using \eqref{eq:XNormOnDeltaiEst} and again the size of $\Delta(t_i)$, we obtain
from \eqref{eq:XNormDiv}
\begin{equation*}
\begin{split}
 \int_0^\infty\|T_{\sigma_i(\cdot)}(\tau)x\|_X^p\,d\tau	&\geq \int_{\Delta(t_i)} \|T_{\sigma_i(\cdot)}(\tau)x\|_X^p\,d\tau\\
						&\geq \left(\frac{\delta_i}{Me}\right)^p\frac{1}{\omega} \to \infty
\end{split}
\end{equation*}
as $i \to \infty$. This contradicts assumption ii). Hence, \eqref{eq:XNormUBound} holds true.

Next, suppose \eqref{eq:StrongStab} does not hold. Then, there exist $x\in X$, 
$\sigma(\cdot)\in\Sigma$, $\delta>0$ and a diverging sequence of times $(t_i)_{i\in\NN}$ such that
\begin{equation}\label{eq:NoStrongStab}
 \|T_{\sigma(\cdot)}(t_i)x\|_X \geq \delta \quad\text{for all}~i.
\end{equation}
Without loss of generality $t_i>t_{i-1}+\frac1\omega$ for every $i\in\mathbb{N}$ 
with $t_0=0$. 
For $\tau\in\Delta(t_i)$, 
\eqref{eq:XY} yields
\begin{equation}\label{eq:XNormOnDeltaiEst2}
 \|T_{\sigma(\cdot)}(\tau)x\|_X \geq \frac{\delta}{Me}, \quad\tau\in\Delta(t_i).
\end{equation}
Hence, using \eqref{eq:XNormOnDeltaiEst} and again the size of $\Delta(t_i)$, we obtain
\begin{equation*}
\begin{split}
 \int_0^\infty\|T_{\sigma(\cdot)}(\tau)x\|_X^p\,d\tau	&\geq \sum_{i=1}^\infty \int_{\Delta(t_i)} \|T_{\sigma(\cdot)}(\tau)x\|_X^p\,d\tau\\
						&\geq  \left(\frac{\delta}{Me}\right)^p\sum_{i=1}^\infty\frac{1}{\omega} =  \infty.
\end{split}
\end{equation*}
This again contradicts assumption ii) and hence \eqref{eq:StrongStab} holds true.

Let
\begin{equation*}
 t_{x,\sigma(\cdot)}(\rho)=\max\{ t : \|T_{\sigma(\cdot)}(t)x\|_X \geq \rho \|x\|_X,~0 \leq s \leq t\}.
\end{equation*}
By \eqref{eq:StrongStab}, $t_{x,\sigma(\cdot)}(\rho)$ is finite (and positive) for every $\sigma(\cdot)$ and $x \in X\setminus\{0\}$.
By strong continuity,
\begin{equation*}
 \|T_{\sigma(\cdot)}(t_{x,\sigma(\cdot)}(\rho))x\|_X = \rho \|x\|_X.
\end{equation*}
Using assumption ii),
\begin{equation*}
\begin{split}
 t_{x,\sigma(\cdot)}(\rho)\rho^p\|x\|_X^p 	&\leq \int_0^{t_{x,\sigma(\cdot)}(\rho)}\|T_{\sigma(\cdot)}(t)x\|_X^p\,dt\\
						&\leq \int_0^{\infty}\|T_{\sigma(\cdot)}(t)x\|_X^p\,dt\leq C\|x\|_X^p
\end{split}
\end{equation*}
whereby,
\begin{equation*}
 t_{x,\sigma(\cdot)}(\rho) \leq \frac{C}{\rho^p}=:t_0,~\text{independent of}~\sigma(\cdot).
\end{equation*}

By the principle of uniform boundedness, \eqref{eq:XNormUBound} implies that there exists 
a constant $k>0$ 
such that
\begin{equation}\label{eq:kdef}
 \|T_{\sigma(\cdot)}(t)\|_{\Lcal(X)} \leq k, \quad t \geq 0,~\sigma(\cdot)\text{-uniformly}.
\end{equation}
Hence, for $t>t_0$, we have
\begin{align*}
 \|T_{\sigma(\cdot)}(t)x\|_X &\leq \sup_{\tilde\sigma(\cdot)}\|T_{\tilde\sigma(\cdot)}(t-t_{x,\sigma(\cdot)}(\rho))\|_{\Lcal(X)}\rho\|x\|_X\\
 &\leq k\rho\|x\|_X, \quad \sigma(\cdot)\text{-uniformly}.
\end{align*}
Choose $\rho>0$ such that $\beta:=k\rho<1$, so that 
\begin{equation*}
 \|T_{\sigma(\cdot)}(t)x\|_X \leq \beta \|x\|_X, \quad t \geq t_0,~\sigma(\cdot)\text{-uniformly}.
\end{equation*}
Finally, let $t_1>t_0$ be fixed and let $t=nt_1+s$, $0 \leq s < t_1$. Then,
 \begin{align*}
  \|T_{\sigma(\cdot)}(t)\|_{\Lcal(X)} &\leq \sup_{\tilde\sigma(\cdot),\;\hat \sigma(\cdot)} \|T_{\tilde\sigma(\cdot)}(s)\|_{\Lcal(X)} \|T_{\hat\sigma(\cdot)}(nt_1)\|_{\Lcal(X)}\\
& \leq k \left(\sup_{\hat\sigma(\cdot)} \|T_{\hat\sigma(\cdot)}(t_1)\|_{\Lcal(X)}\right)^n\\
& \leq k \beta^n \leq K e^{-\mu t},~t \geq 0,~\sigma(\cdot)\text{-uniformly}
 \end{align*}
with $K=\frac{C}{\beta}$ and $\mu=-\left(\frac{1}{t_1}\right)\ln\beta>0$.
\end{proof}

Lemma~\ref{lem:IntUBound} allows us to prove the first of the converse Lyapunov theorems that are the object of this paper.  
\begin{theorem}\label{thm:VnonComparable} 
The following conditions:
\begin{itemize}
\item[i)] there exist constants $M \geq 1$ and $\omega>0$
such that
\begin{equation}\label{thm1:UEBound}
 \|T_{\sigma(\cdot)}(t)\|_{\Lcal(X)} \leq Me^{\omega t}, \quad t \geq 0,~\sigma(\cdot)\text{-uniformly},
\end{equation}
 \item[ii)] 
 there exists $V\: X \to [0,\infty)$
 such that 
 $\sqrt{V(\cdot)}$ is a norm on $X$, 
	\begin{equation}\label{thm1:VBound}
 		V(x) \leq C \|x\|_X^2, \quad x\in X,
	\end{equation}
	for a constant $C>0$ and
	\begin{equation}\label{thm1:LjVest}
 		\Linf_j V(x) \leq -\|x\|_X^2,\quad j \in Q,~x\in X,
	\end{equation}
	with $\Linf_j V(x)$ defined as in \eqref{eq:LinfjVdef}
\end{itemize}
are necessary and sufficient for the existence of constants $K \geq 1$ and $\mu>0$
such that
\begin{equation}\label{thm1:GUES}
 \|T_{\sigma(\cdot)}(t)\|_{\Lcal(X)} \leq Ke^{-\mu t}, \quad t \geq 0,~\sigma(\cdot)\text{-uniformly}.
\end{equation}
\mbox{}
\end{theorem}
\begin{proof}
Assume that the conditions i) and ii) hold. For all $\sigma(\cdot)\in\Sigma$, $x \in X$
and for $t \geq 0$ small enough so that the restriction of $\sigma(\cdot)$ to the interval $[0,t]$ is constant, we have
\begin{equation*}
 0 \leq V(T_{\sigma(\cdot)}(t)x) \leq V(x) - \int_0^t \|T_{\sigma(\cdot)}(\tau)x\|_X^2\,d\tau
\end{equation*}
as it follows from \eqref{thm1:LjVest} and \cite[\S VI.7]{Saks1964} (see also \cite{HagoodThomson2006}). Thus, for all $\sigma(\cdot)$ and $x \in X$,
\begin{equation}\label{eq:IntUBound}
\int_0^\infty \|T_{\sigma(\cdot)}(\tau)x\|_X^2\,d\tau \leq V(x)\leq C\|x\|_X^2.
\end{equation}
The uniform exponential decay \eqref{thm1:GUES} now follows from \eqref{thm1:UEBound} and \eqref{eq:IntUBound}, thanks to
Lemma~\ref{lem:IntUBound} with $p=2$.

Conversely, assume that \eqref{thm1:GUES} holds for some constants $K \geq 1$ and $\mu>0$.
Then, \eqref{thm1:UEBound} holds for $M=K$ and arbitrary $\omega>0$. Define $V\: X \to [0,\infty)$ by
\begin{equation}\label{eq:Vdef}
V(x) = \sup_{\sigma(\cdot)\in\Sigma}\int_0^\infty \|T_{\sigma(\cdot)}(t)x\|_X^2\,dt. 
\end{equation}
Then, by assumption, $V(x)$ satisfies
\begin{equation*}
 V(x) \leq \sup_{\sigma(\cdot)}\int_0^\infty K^2e^{-2\mu t}\|x\|_X^2\,dt=\frac{K^2}{2\mu}\|x\|_X^2,
\end{equation*}
establishing \eqref{thm1:VBound} with $C=\frac{K^2}{2\mu}$. In particular, $V$ is well-posed.

Notice that, by definition, $V$ is positive definite and homogenous of degree $2$.
In order to show that it is the square of a norm, we are left to prove that it is convex and continuous.

The convexity of $V$ follows from the fact that each 
$$x\mapsto \int_0^\infty \|T_{\sigma(\cdot)}(t)x\|_X^2\,dt$$
is convex.

In order to verify the continuity of $V$, let $(x_n)_{n\in\NN}$ be a sequence in $X$ converging to $x$ in $X$. By definition of $V$,
\begin{equation}\label{eq:Vnest}
\int_0^\infty \|T_{\sigma(\cdot)}(t)x_n\|_X^2\,dt \leq V(x_n),
\end{equation}
for all $\sigma(\cdot)\in\Sigma$. So taking the $\liminf$ over $n\in\NN$ in \eqref{eq:Vnest} on both sides 
and using the continuity of $T_{\sigma(\cdot)}(t)$ for all $t\geq 0$, we have
\begin{equation}\label{eq:Vcont}
 \int_0^\infty \|T_{\sigma(\cdot)}(t)x\|_X^2\,dt \leq \liminf_{n\to\infty} \int_0^\infty \|T_{\sigma(\cdot)}(t)x_n\|_X^2\,dt \leq \liminf_{n\to\infty} V(x_n).
\end{equation}
Taking the $\sup$ over $\sigma(\cdot)$ in \eqref{eq:Vcont} then yields
\begin{equation*}
 V(x) = \sup_{\sigma(\cdot)}\int_0^\infty \|T_{\sigma(\cdot)}(t)x\|_X^2\,dt \leq \liminf_{n\to\infty} V(x_n),
\end{equation*}
proving that $V$ is lower semi-continuous. On the other hand, for a fixed $\varepsilon>0$, there exist $\sigma_{\varepsilon}(\cdot)\in\Sigma$ such that
\begin{align*}
 V(x_n) - \frac{\varepsilon}{2} &< \int_0^\infty\kern-0.4em\|T_{\sigma_{\varepsilon}(\cdot)}(t)x_n\|^2\,dt\\
&\leq (1+m)\int_0^\infty\kern-0.4em \|T_{\sigma_{\varepsilon}(\cdot)}(t)(x_n-x)\|^2\,dt+ \left(1+\frac1m\right)\int_0^\infty\kern-0.4em \|T_{\sigma_{\varepsilon}(\cdot)}(t)x\|^2\,dt,
\end{align*}
for any $m>0$. 
Notice that, by definition of $V$,
\begin{align*}
 \int_0^\infty \|T_{\sigma_{\varepsilon}(\cdot)}(t)(x_n-x)\|^2\,dt&\leq
 V(x_n-x)
 \leq C\|x_n-x\|_X^2,\\
 \int_0^\infty \|T_{\sigma_{\varepsilon}(\cdot)}(t)x\|^2\,dt &\leq V(x).
\end{align*}
Thus, for any $m>0$, we have
\begin{equation}\label{eq:uscestLast}
V(x_n) - \frac{\varepsilon}{2} < C (1+m)\|x_n-x\|^2_X + \left(1+\frac1m\right)V(x).
\end{equation}
In particular, choosing $m$ such that $(1+\frac1m)V(x)<V(x)+\frac\varepsilon4$ and taking $n$ sufficiently 
large, so that $(1+m) C\|x_n-x\|_X^2 \leq \frac\varepsilon4$, we have from \eqref{eq:uscestLast}
\begin{equation*}
 V(x_n) < V(x) + \varepsilon, \quad n~\text{sufficienty large}.
\end{equation*}
This implies the upper semi-continuity of $V$. 
Resuming, we proved the
continuity of $V$. 

To complete the proof of the theorem, we are left to show that $V$ satisfies \eqref{thm1:LjVest}. 
Fixing $t > 0$, $j \in Q$ and letting
\begin{equation*}
\Sigma_{t,j}=\{\sigma(\cdot)\in\Sigma : \sigma|_{[0,t]}\equiv j\}
\end{equation*}
be the set of switching signals whose restriction to the interval $[0,t]$ is constantly equal to $j$, we have, since $\Sigma_{t,j}\subseteq\Sigma$,
\begin{align}
 V(x)	&\geq\sup_{\sigma(\cdot)\in\Sigma_{t,j}} \int_0^\infty \|T_{\sigma(\cdot)}(\tau)x\|_X^2\,d\tau\nonumber\\
 	&=\int_0^t \|T_{j}(\tau)x\|_X^2\,d\tau + \sup_{\sigma(\cdot)\in\Sigma_{t,j}} \int_t^\infty \|T_{\sigma(\cdot)}(\tau)x\|_X^2\,d\tau.\label{half}
\end{align}
Moreover, thanks to \eqref{eq:SemigoupProperty} and the invariance of $\Sigma$ by time-shift,
\begin{align*}
V(T_j(t)x)	&=\sup_{\sigma(\cdot)\in\Sigma} \int_0^\infty \|T_{\sigma(\cdot)}(\tau)T_j(t)x\|_X^2\,d\tau\\
		&=\sup_{\sigma(\cdot)\in\Sigma} \int_t^\infty \|T_{\sigma(\cdot)}(\tau-t)T_j(t)x\|_X^2\,d\tau\\
 		&=\sup_{\sigma(\cdot)\in\Sigma_{t,j}} \int_t^\infty \|T_{\sigma(\cdot)}(\tau)x\|_X^2\,d\tau.
\end{align*}
This and \eqref{half} yield
\begin{equation*}
 V(T_j(t)x)-V(x) \leq - \int_0^t \|T_j(\tau)\|_X^2\,d\tau,
\end{equation*}
for all $j\in Q$ and $t > 0$. Therefore
\begin{equation*}
\begin{split}
 &\Linf_j V(x) = \liminf_{t \downarrow 0} \frac{V(T_j(t)x) - V(x)}{t}\\
 &\qquad\leq - \limsup_{t \downarrow 0} \frac{1}{t}\int_0^t \|T_j(\tau)x\|_X^2\,d\tau=-\|x\|_X^2
\end{split}
\end{equation*}
for all $j\in Q$, establishing \eqref{thm1:LjVest}.
\end{proof}

\begin{figure*}[t]
\centering
\scalebox{1}[1]{
\begin{pspicture}(5,5)
\pspolygon[linestyle=none,fillstyle=crosshatch,hatchangle=45,hatchwidth=0.1pt](1,4)(1,1.5)(1.5,1)(1.5,3.5)(1,4)
\pspolygon[linestyle=none,fillstyle=vlines,hatchangle=45,hatchwidth=0.1pt](1,1)(1.5,1)(1,1.5)(1,1)
\pspolygon[linestyle=none,fillstyle=vlines,hatchangle=45,hatchwidth=0.1pt](1.5,1)(4,1)(1.5,3.5)(1.5,1)
\psline{->}(1,1)(4.5,1)
\psline{->}(1,1)(1,4.5)
\psline[linestyle=dotted]{-}(4,1)(4,4.5)
\psline[arrowsize=0.2,linewidth=0.5pt]{->}(4,1)(1,4)
\psline[arrowsize=0.2,linewidth=0.5pt]{->}(1.5,1)(1,1.5)
\psline[linestyle=dotted]{-}(1.5,1)(1.5,4.5)
\rput(0.5,2.5){$t$}
\rput(5,0.75){$s$}
\psline[linestyle=dotted]{-}(0.95,4)(4,4)
\psline{-}(0.95,4)(1.05,4)
\rput(0.75,4){$1$}
\psline{-}(1,0.95)(1,1.05)
\rput(0.9,0.75){$0$}
\psline{-}(1.5,0.95)(1.5,1.05)
\rput(1.5,0.775){$4^{-j}$}
\psline{-}(4,0.95)(4,1.05)
\rput(4,0.75){$1$}
\end{pspicture}
}
\caption{Illustration of $T_j(t)$, $j \in Q$, of the example in Remark~\ref{rem:Thm1sharp}.\label{f::2}}
\end{figure*}
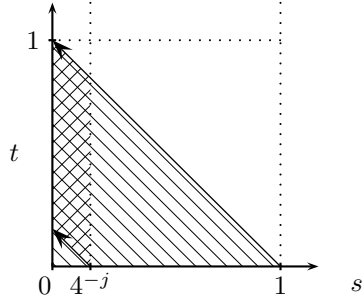

\begin{remark}\label{rem:Thm1sharp} We show here, through an example, that  condition i) appearing in the statement of Theorem~\ref{thm:VnonComparable} 
cannot be removed.

Consider the family of semigroups $\{T_j(t)\}_{j \in Q}$ with $Q=\NN$ and $T_j(t)$ defined on the Lebesgue space
$X=L^p(0,1)$, $p\in[1,\infty)$, by
\begin{equation}\label{eq:TcounterexDef}
 \left(T_j(t)f\right)(s)=\begin{cases}
                         	2^{\frac1p} f(s+t),	& s \in [0,1-t]\cap[0,{4^{-j}})\\
                         	f(s+t),			& s \in [0,1-t]\setminus[0,4^{-j})\\
				0,			& s \in (1-t,1],
                         \end{cases}
\end{equation}
cf. also Figure~\ref{f::2}. Notice that, for all $j \in Q$, $T_j(\cdot)$ is a nilpotent semigroup, since $\|T_j(t) f\|_X=0$ for $t > 1$. In particular, each of them is exponentially stable. Moreover, for all $\sigma(\cdot)\in\Sigma$, we have 
\begin{equation*}
\left|\left(T_{\sigma(\cdot)}(t)f\right)(s)\right|\leq 2^{\frac k p}|f(s+t)|,\quad\mbox{with }k=\#\{l\in\NN :
 s<4^{-l}\leq s+t\}.
\end{equation*}
 In particular,
 \begin{equation*}
\left|\left(T_{\sigma(\cdot)}(t)f\right)(s)\right|\leq 2^{\frac k p}|f(s+t)|,\quad\mbox{if }  4^{-k-1}\leq s< 4^{-k},
\end{equation*}
 yielding
\begin{align*}
 \int_0^\infty \|T_{\sigma(\cdot)}(t)f\|_X^p 
 &= \int_0^1 \int_0^{1-t}\left|\left(T_{\sigma(\cdot)}(t)f\right)(s)\right|^p\,ds\,dt\\
&= \int_0^1 \int_0^{1-s}\left|\left(T_{\sigma(\cdot)}(t)f\right)(s)\right|^p\,dt\,ds\\
&\leq \sum_{k=0}^\infty \int_{4^{-k-1}}^{4^{-k}} 2^{k}\left( \int_0^{1-s}| f(s+t)|^p\,dt\right)ds\\
&\leq
\sum_{k=0}^\infty \left(\frac{1}{{4^{k}}}-\frac{1}{{4^{k+1}}}\right)2^{k}\int_0^1|f(t)|^p\,dt = \frac{3}{2}\|f\|_X^p.
\end{align*}

Hence, defining $V(x)$ as in \eqref{eq:Vdef}, we have
\begin{equation*}
 V(x)=\sup_{\sigma(\cdot)} \int_0^\infty \|T_{\sigma(\cdot)}(t)x\|_X^2\,dt \leq \frac{3}{2}\|x\|_X^2
\end{equation*}
and, by the same arguments as in the proof of Theorem~\ref{thm:VnonComparable}, 
\begin{equation*}
\Linf_j V(x) \leq -\|x\|_X^2, \quad\text{for all}~j \in Q, 
\end{equation*}
so condition ii) in Theorem~\ref{thm:VnonComparable} holds with $C=\frac{3}{2}$.

Nevertheless, for a sequence of switching signals $(\sigma_n(\cdot))_{n\in\NN}\subset\Sigma$ with switching times $\tau_k=\frac{1}{4^k}$ and modes $j_k=k+1$, $0\leq k\leq n$, we have for functions $\One_{[s,1]}$ of $L^p$ norm one 
concentrated on the interval $[s,1]$,
\begin{equation*}
 T_{\sigma_n(\cdot)}(1-\epsilon) \One_{[s,1]} = 2^{\frac{n}{p}}
  \One_{[s-1+\epsilon,\epsilon]}, \quad\mbox{if }1\geq s\geq 1-\epsilon>1-4^{-n}.
\end{equation*}
Therefore, for $\epsilon<4^{-n}$, 
\begin{align*}
\|T_{\sigma_n(\cdot)}(1-\epsilon)\|_{\Lcal(X)} 	&= \sup_{\|f\|_X=1} \|T_{\sigma_n(\cdot)}(1-\epsilon)f\|_X\\
						&\geq \lim_{s \uparrow 1} \|T_{\sigma_n(\cdot)}(1-\epsilon) \One_{[s,1]}\|_X
						= 2^{\frac n p}.
\end{align*}

Hence, 
\begin{equation}\label{eq:exampleblowup}
\sup\{\|T_{\sigma_n(\cdot)}(1-\epsilon)\|_{\Lcal(X)} : \epsilon\in[0,1],~n\in\NN\}=+\infty 
\end{equation}
violating any uniform bound of the form \eqref{thm1:GUES}. 

This example also shows that assumption i) appearing in Lemma~\ref{lem:IntUBound} is 
necessary for the validity of the lemma. 
\hfill$\square$
\end{remark}

\section{SECOND CONVERSE LYAPUNOV THEOREM}\label{sec:second}

If one wishes to conclude that a switched system is globally uniformly exponentially stable, Theorem~\ref{thm:VnonComparable}
requires the knowledge of a squared norm $V(\cdot)$ satisfying \eqref{thm1:VBound},
\eqref{thm1:LjVest} for all modes $j \in Q$ and the knowledge of a global uniform exponential bound \eqref{thm1:UEBound}. 
As an alternative, we will, in Theorem~\ref{thm:VComparable}, show that the 
existence of a Lyapunov norm $\sqrt{V(\cdot)}$ that is comparable with the norm $\|\cdot\|_X$ allows to 
conclude that the system is globally uniformly exponentially stable, without knowledge of a global uniform exponential 
bound of the type \eqref{thm1:UEBound}.
Notice that
the norm $\sqrt{V(\cdot)}$ constructed in the proof of Theorem~\ref{thm:VnonComparable} (see definition \eqref{eq:Vdef}) is in general not comparable with 
$\|\cdot\|_{X}$, i.\,e., in general 
it
does not satisfy a lower bound of the form
\begin{equation}\label{lower_est}
 c\|x\|_X \leq \sqrt{V({x})}, \quad x \in X
\end{equation}
for a constant $c>0$. 
Such a lower bound always holds, on the contrary, 
when $X$ has finite dimension,
as it is exploited in \cite{MolchanovPyatnitskiy1986,DayawansaMartin1999}. 
The bound \eqref{lower_est} may fail to hold even in the case of a single strongly continuous semigroup, as it is the case, for instance, of the semigroups $T_1(\cdot)$ and $T_2(\cdot)$ introduced in Example~\ref{ex:notGUES}.   (For a characterization of 
exponentially stable strongly continuous semigroups 
whose Lyapunov function defined as in \eqref{eq:Vdef}
is comparable with the squared norm, see \cite{Pazy1972}.)

In order to obtain a Lyapunov norm 
comparable with 
$\|\cdot\|_X$
for infinite dimensional switched 
systems, we make use of the following lemma imposing a stronger assumption on the
family of semigroups $T_j(\cdot)$.

\begin{lemma}\label{lem:VComparable}
Assume that there exists $j^\ast \in Q$ such that $T_{j^\ast}(\cdot)$ can be extended to a group of bounded 
linear operators on $X$. Moreover, assume that there exist constants $K \geq 1$ and $\mu>0$ 
such that
\begin{equation}\label{lem2:GUES}
 \|T_{\sigma(\cdot)}(t)\|_{\Lcal(X)} \leq Ke^{-\mu t}, \quad t \geq 0,~\sigma(\cdot)\text{-uniformly}.
\end{equation}
Then there exists 
a 
function 
$V\: X \to [0,\infty)$ 
such that 
$\sqrt{V(\cdot)}$ is a norm on $X$,
\begin{equation}\label{lem2:VBounds}
	c\|x\|_X^2 \leq V(x) \leq C \|x\|_X^2,\quad x\in X,
\end{equation}
for constants $c,C>0$ and
\begin{equation}\label{lem2:LjVest}
	\Linf_j V(x) \leq -\|x\|_X^2, \quad j \in Q,~x\in X,
\end{equation}
with $\Linf_j V(x)$ defined as in \eqref{eq:LinfjVdef}.
\end{lemma}
\begin{proof}
Assume that \eqref{lem2:GUES} holds for some constants $K \geq 1$ and $\mu>0$ independent of $\sigma(\cdot)$.
Define $V(\cdot)$ by \eqref{eq:Vdef}. As seen in 
the proof of Theorem~\ref{thm:VnonComparable}, 
\eqref{lem2:GUES} guarantees that 
$V$ is 
the square of a norm,
satisfies \eqref{lem2:LjVest} and that there exist $C>0$ such that 
$V(x) \leq C \|x\|_X^2$. 

The remaining bound $c\|x\|^2 \leq V(x)$ for some constant $c>0$ is a consequence of the assumption that $T_{j^\ast}(\cdot)$  can be extended to a group. Indeed, $T_{j^\ast}(t)$ is then invertible for every $t \geq 0$ and satisfies $\|T_{j^\ast}(t)x\|_X \geq (\|T_{j^\ast}(-t)\|_{\Lcal(X)})^{-1}\|x\|_X$  (cf.~\cite{Pazy1972}). Hence
\begin{equation*}
\begin{split}
 V(x)	&= \sup_{\sigma(\cdot)} \int_0^\infty \|T_{\sigma(\cdot)}(t)x\|_X^2\,dt \geq \int_0^\infty \|T_{j^\ast}(t)x\|_X^2\\
	&\geq \int_0^\infty (\|T_{j^\ast}(-t)\|_{\Lcal(X)})^{-2}\,dt \|x\|_X^2 = c\|x\|_X^2
\end{split}
\end{equation*}
with $c=\int_0^\infty (\|T_{j^\ast}(-t)\|_{\Lcal(X)})^{-2}\,dt<\infty$. 
\end{proof}

We can now state and prove our second converse Lyapunov theorem. 

\begin{theorem}\label{thm:VComparable} The existence of  
a 
function 
$V\: X \to [0,\infty)$ 
such that
$\sqrt{V(\cdot)}$ is a norm on $X$,
\begin{equation}\label{thm2:VBounds}
	c\|x\|_X^2 \leq V(x) \leq C \|x\|_X^2, \quad x\in X,
\end{equation}
for constants $c,C>0$ and
\begin{equation}\label{thm2:LjVest}
	\Linf_j V(x) \leq -\|x\|_X^2, \quad j \in Q,~x\in X,
\end{equation}
with $\Linf_j V(x)$ defined as in \eqref{eq:LinfjVdef} is necessary and sufficient for
the existence of constants $K \geq 1$ and $\mu>0$
such that
\begin{equation}\label{thm2:GUES}
 \|T_{\sigma(\cdot)}(t)\|_{\Lcal(X)} \leq Ke^{-\mu t}, \quad t \geq 0,~\sigma(\cdot)\text{-uniformly}.
\end{equation}
\mbox{}
\end{theorem}
\begin{proof}
Assume that there exists a function $V\: X \to [0,\infty)$ such that \eqref{thm2:VBounds}
and \eqref{thm2:LjVest} hold for 
$c,C>0$. Then, for all $x \in X$, $\sigma(\cdot)\in\Sigma$, and $t \geq 0$,
\begin{equation*}
 c\|T_{\sigma(\cdot)}(t)x\|_X^2 \leq V(T_{\sigma(\cdot)}(t)x) \leq V(x) - \int_0^t \|T_{\sigma(\cdot)}(\tau)x\|_X^2\,d\tau,
\end{equation*}
so, using again \eqref{thm2:VBounds} and dividing by $c$,
\begin{equation*}
 \|T_{\sigma(\cdot)}(t)x\|_X^2 \leq \frac{C}{c}\|x\|_X^2 - \frac{1}{c}\int_0^t \|T_{\sigma(\cdot)}(\tau)x\|_X^2\,d\tau.
\end{equation*}
From Gronwall's Lemma, we obtain
\begin{equation*}
 \|T_{\sigma(\cdot)}(t)x\|_X^2 \leq \frac{C}{c} e^{-\frac{1}{c}t} \|x\|_X^2, \quad t \geq 0,~\sigma(\cdot)\text{-uniformly}.
\end{equation*}
Hence, 
\begin{equation*} 
\|T_{\sigma(\cdot)}(t)\|_{\Lcal(X)} \leq Ke^{-\mu t}, \quad t \geq 0,~\sigma(\cdot)\text{-uniformly},
\end{equation*}
for the constants $K=\sqrt{\frac{C}{c}} \geq 1$ and $\mu=\frac{1}{2c}>0$, establishing \eqref{thm2:GUES}.

Conversely, assume that \eqref{thm2:GUES} holds for 
$K \geq 1$ and $\mu>0$ and consider the 
switched system with $\sigma^\ast(\cdot)$ taking values in $Q^\ast=Q\cup\{j^*\}$, where
\begin{equation}\label{def:SSextended}
 T_{j^\ast}(t)=e^{-\mu t} I
\end{equation}
and $I\: X \to X$ denotes the identity.
Then, 
\begin{equation*}
 T_{j^\ast}(t) T_j(s) = T_j(s) T_{j^\ast}(t), \quad t,s \geq 0,~j \in Q.
\end{equation*}
Moreover, for $t^\ast=|\{\tau \in [0,t] : \sigma^\ast(\tau)=j^\ast\}|$ and
some $\sigma(\cdot)$ just taking values in $Q$,
\begin{align}\label{eq:SSextendedGUES}
\|T_{\sigma^\ast(\cdot)}\|_{\Lcal(X)} 
&= \|T_{\sigma(\cdot)}(t-t^\ast)e^{-\mu t^\ast}\|_{\Lcal(X)} \leq K e^{-\mu(t-t^\ast)} e^{-\mu t^\ast}\nonumber\\
					&=K e^{-\mu t}, \quad\quad t \geq 0,~\sigma(\cdot)\text{-uniformly}
\end{align}
where we have used \eqref{thm2:GUES}.
The existence of a squared norm 
$V\: X \to [0,\infty)$ 
and constants $c,C>0$
such that \eqref{thm2:VBounds}
and \eqref{thm2:LjVest} hold 
now follows from \eqref{eq:SSextendedGUES}
and Lemma~\ref{lem:VComparable}, noting that 
\eqref{def:SSextended}
actually defines 
 a group.
\end{proof}

\begin{remark}\label{one-mode-L}
In the case of a single strongly continuous semigroup $T(\cdot)$, 
Theorem~\ref{thm:VComparable} shows
that its global exponential stability is equivalent to 
the existence of a Lyapunov norm comparable with the norm $\|\cdot\|_{X}$.  

When $T(\cdot)$ is globally exponentially stable, such a Lyapunov norm can be obtained 
by 
the construction suggested in the proof of Theorem~\ref{thm:VComparable}.

An alternative construction of a common Lyapunov function has been proposed in \cite{HanteSigalotti2010}, where
\eqref{eq:Vdef}
is replaced by 
\begin{equation}\label{defV:CDC}
\tilde V(x) = \int_0^\infty \sup_{\sigma(\cdot)\in\Sigma}\|T_{\sigma(\cdot)}(t)x\|_X^2\,dt. 
\end{equation}
In the case of a single mode and following the strategy of augmenting $Q$ by adding a group corresponding to  a diagonal operator, this construction leads to the explicit expression
\begin{equation*}
\tilde V(x)=\int_0^\infty \left(\max_{s\in[0,\tau]}e^{2\mu(s-\tau)}\|T(s)x\|_X^2\right)d\tau,
\end{equation*}
for any fixed $\mu>0$. 

It should be noticed that, although not stated in \cite{HanteSigalotti2010}, the definition of $\tilde V$ given in \eqref{defV:CDC} 
identifies a function which is positive definite, homogenous of degree 2,   
continuous and convex, i.e., a squared norm. The proof 
of this fact can be rather easily obtained by 
adapting the proof of Theorem~\ref{thm:VnonComparable}.
\hfill$\square$
\end{remark}

\section{COMMON LYAPUNOV FUNCTIONS ON 
HILBERT SPACES}\label{sec:Hilbert}
The goal of this section is to prove further regulatity properties 
of the common Lyapunov functions constructed in the previous sections 
in the case in which $X$ is a separable Hilbert space. 
The proof of the following lemma adapts 
arguments presented in \cite[\S4.3.1]{BonnansShapiro2000}. 

Recall that a function $V\: X\to \RR$ is said to be {\it directionally differentiable in the sense of Fr\'{e}chet}
at $x\in X$ if for every $\psi\in X$ there exists
$$V'(x,\psi)=\lim_{t\downarrow 0}\frac{V(x+t\psi)-V(x)}{t}$$
and, moreover,
$$\lim_{\psi\to 0}\frac{V(x+\psi)-V(x)-V'(x,\psi)}{\|\psi\|_X}=0.$$
Notice that, in general, $V'(x,\cdot)$ needs not be linear.

\begin{lemma} \label{lem:Hilbert} 
Let $X$ be a separable Hilbert space with scalar product $\langle\cdot,\cdot\rangle$ and assume that there exist constants $K \geq 1$ and $\mu>0$
such that
\begin{equation}\label{prop1:GUES}
 \|T_{\sigma(\cdot)}(t)\|_{\Lcal(X)} \leq Ke^{-\mu t}, \quad t \geq 0,~\sigma(\cdot)\text{-uniformly}.
\end{equation}
Then, there exists a  
subset $\B$ of $\Lcal(X)$, compact
for the weak operator topology and 
made 
of self-adjoint operators, such  that 
\begin{equation}\label{defVviaB}
 \sup_{\sigma(\cdot)\in\Sigma}\int_0^\infty \|T_{\sigma(\cdot)}(t)x\|_X^2\,dt =\max_{B\in\B} \la x, B x \ra=:V(x).
\end{equation}
In particular, $V$ is 
directionally differentiable in the sense of Fr\'{e}chet and its derivative in the direction $\psi \in X$ is given by
\begin{equation}\label{eq:VDdef}
  V'(x,\psi) = \max_{\hat{B}\in\Scal(x)} 2\la \psi,\hat{B} x\ra,
\end{equation}
where 
\begin{equation}\label{eq:Sdef}
 \Scal(x)=\argmax_{B\in\B} \la x, B x \ra.
\end{equation}
\end{lemma}
\begin{proof}
For a fixed $\sigma(\cdot)\in\Sigma$ and all $t\geq 0$, let $T_{\sigma(\cdot)}^*(t)\in\Lcal(X)$ be the adjoint operator of $T_{\sigma(\cdot)}(t)\in\Lcal(X)$, uniquely defined by
\begin{equation*}
 \langle T_{\sigma(\cdot)}(t)x,x'\rangle = \langle x,T_{\sigma(\cdot)}^*(t)x'\rangle, \quad x,x'\in X.
\end{equation*}
If $T_{\sigma(\cdot)}(t)$ has the expression given in \eqref{eq:TsigmaDef},
then
\begin{equation*}
 T_{\sigma(\cdot)}^*(t)=T_{j_{1}}^*(\tau_{1})\cdots T_{j_{p-1}}^*(\tau_p-\tau_{p-1})T_{j_p}^*(t-\tau_p),
\end{equation*}
where $T_{j}^*(t)$ denotes the adjoint semigroup of $T_{j}(t)$, $t \geq 0$, $j \in Q$. 

Assuming that there exist constants $K \geq 1$ and $\mu>0$, independent of $\sigma(\cdot)$, such that
\eqref{prop1:GUES} holds, the operator $B_{\sigma(\cdot)}\: X \to X$, given by
\begin{equation*}
 B_{\sigma(\cdot)}x = \int_0^\infty T_{\sigma(\cdot)}^*(t)T_{\sigma(\cdot)}(t)x\,dt,
\end{equation*}
is linear, self-adjoint, and satisfies
\begin{align*}
 \|B_{\sigma(\cdot)}x\|_X	&\leq \int_0^\infty \|T_{\sigma(\cdot)}^*(t)\|_{\Lcal(X)}\|T_{\sigma(\cdot)}(t)\|_{\Lcal(X)}\|x\|_X\,dt\\
				&\leq \int_0^\infty K^2e^{-2\mu t}\|x\|_X\,dt=\frac{K^2}{2\mu}\|x\|_X,
\end{align*}
where we have used \eqref{prop1:GUES}. In particular, $B_{\sigma(\cdot)}\in\Lcal(X)$ for all $\sigma(\cdot)\in\Sigma$ and  
\begin{equation}\label{eq:BsigmanBound}
 \|B_{\sigma(\cdot)}\|_{\Lcal(X)} \leq \frac{K^2}{2\mu}, \quad\sigma(\cdot)\text{-uniformly}.
\end{equation}
Therefore, the set
\begin{equation*}
 \B=\{ B \in \Lcal(X) :~\text{there exists a sequence}~(\sigma_n(\cdot))_{n\in\NN}\subset \Sigma~\text{such that}~ B_{\sigma_n(\cdot)} \wot B \}
\end{equation*}
is compact for the weak operator topology. We recall that 
$B_{\sigma_n(\cdot)}\wot B$ (i.\,e., $B_{\sigma_n(\cdot)}$ converges to $B$ for the \emph{weak operator topology}) if, for every sequence
$(x_n,y_n)$ converging to $(x,y)$ in $X\times X$, 
we have 
 \begin{equation}\label{wstard}
\lim_{n\to\infty}\langle x_n,B_{\sigma_n(\cdot)}y_n\rangle=\langle x,B y\rangle,
\end{equation}
and that every bounded closed subset of $L(X)$ is sequentially compact for the weak operator topology (see, for instance, \cite[Theorem III.4]{Fabec}). 
Notice that \eqref{wstard} guarantees that $\B$ 
consists of self-adjoint operators.

Define $V$ as in \eqref{defVviaB}. The maximization makes sense because of 
\eqref{wstard} (with $x_n=y_n=x$) and of the compactness of $\B$. 
Moreover,
\begin{align*}
 V(x) &= \max_{B\in\B}\la x,Bx \ra = \sup_{\sigma(\cdot)\in\Sigma}\la x,B_{\sigma(\cdot)}x \ra = \sup_{\sigma(\cdot)\in\Sigma}\la x, \int_0^\infty T_{\sigma(\cdot)}^*(t)T_{\sigma(\cdot)}(t) x\,dt \ra\\
 &=\sup_{\sigma(\cdot)\in\Sigma}\int_0^\infty \la x, T_{\sigma(\cdot)}^*(t)T_{\sigma(\cdot)}(t) x \ra dt 
= \sup_{\sigma(\cdot)}\int_0^\infty \| T_{\sigma(\cdot)}(t) x\|^2\,dt,
\end{align*}
for all $x \in X$. 
Hence $V$ coincides with the Lyapunov function defined in  \eqref{eq:Vdef} and we recover that $V$ satisfies the condition ii) of Theorem~\ref{thm:VnonComparable}. In particular, $V$ is 
continuous. 

In order to verify the directional differentiability in the sense of Fr\'{e}chet and 
to prove \eqref{eq:VDdef}, we show below that
\begin{equation}\label{Frechet_D}
\lim_{\psi\to 0}\frac{V(x+\psi)-V(x)-\max_{\hat{B}\in\Scal(x)} 2\la \psi,\hat{B} x\ra}{\|\psi\|_X}=0.
\end{equation}

First observe that the map
$x \mapsto \la x,Bx \ra$ is differentiable on $X$ in the sense of Fr\'{e}chet for all $B\in\Lcal(X)$. For $B$ self-adjoint, the derivative is given by $2\la \cdot, Bx\ra$. Fix now some $x_0 \in X$ and define
\begin{equation*}
\Phi(B,x) = 2 \|Bx - B x_0\|_{X}.
\end{equation*}
We claim that 
\begin{equation}\label{eq:PhiunifToZero}
 \lim_{x \to x_0} \Phi(B,x)=0, \quad\text{uniformly with respect to}~B\in\B.
\end{equation}
Indeed, suppose by contradiction that there exists some $\varepsilon>0$, a sequence $(x_n)_{n\in\NN}$ converging to $x_0$ in $X$ and a sequence $(B_n)_{n\in\NN}$ in $\B$ such that 
\begin{equation}\label{eq:proofcontra1}
\Phi(B_n,x_n)>\varepsilon \quad\text{for all}~n\in\NN.
\end{equation}
Thanks to the compactness of $\B$, there exist $B\in\B$ and  a subsequence $n(k)$ such that 
\begin{equation*}
B_{n(k)}\wot B\quad\mbox{as $k\to\infty$}. 
\end{equation*}
Then, it follows from  \eqref{wstard} that
\begin{equation*}
 \Phi(B_{n(k)},x_{n(k)}) =  0,
\end{equation*}
contradicting \eqref{eq:proofcontra1}. 

The Mean Value Theorem then gives
\begin{equation}\label{eq:MVTest1}
\begin{split}
&|\la x_0 + \psi, B (x_0 + \psi)\ra - \la x_0, B x_0 \ra - 2 \la \psi,B x_0 \ra| \\
&\qquad\qquad\qquad\qquad\leq \|\psi\|_X \int_0^1 \Phi(B,x_0+\xi\psi)\,d\xi \leq \varepsilon \|\psi\|_X
\end{split}
\end{equation}
for $\varepsilon>0$, $\psi$ sufficiently close to zero, and uniformly with respect to 
$B\in\B$, as it follows from \eqref{eq:PhiunifToZero}.

Let $\Scal(\cdot)$ be defined as in \eqref{eq:Sdef}.
Notice that, for every $x\in X$, $\Scal(x)$ is close and therefore compact for the weak operator topology. 
For any $\hat{B}\in\Scal(x_0)$, 
$V(x_0)=\la x_0,\hat{B}x_0\ra$ and $V(x_0+\psi) \geq \la x_0+\psi,\hat{B}(x_0+\psi)\ra$. Thus, using \eqref{eq:MVTest1},
\begin{equation*}
 V(x_0+\psi) - V(x_0) \geq \max_{\hat{B}\in\Scal(x_0)} 2\la \psi,\hat{B} x_0\ra - \varepsilon \|\psi\|_X.
\end{equation*}
In order to prove \eqref{Frechet_D}, it therefore remains to show that for all $\varepsilon>0$ and $\psi$ close to zero,
\begin{equation}\label{eq:proofcontra2}
 V(x_0+\psi) - V(x_0) \leq \max_{\hat{B}\in\Scal(x_0)} 2\la \psi,\hat{B} x_0\ra + \varepsilon \|\psi\|_X.
\end{equation}
So, suppose \eqref{eq:proofcontra2} not to hold. Then there exist $\varepsilon>0$ and a sequence $(\psi_n)_{n\in\NN}$ in $X$ converging to zero such that
\begin{equation}\label{eq:proofcontra3}
 V(x_0+\psi_n) - V(x_0) > \max_{\hat{B}\in\Scal(x_0)} 2\la \psi_n,\hat{B} x_0\ra + \varepsilon \|\psi_n\|_X.
\end{equation}
By definition of $V$, 
there exist $B_0,\hat{B}_n\in\B$
such that
\begin{align*}
 V(x_0+\psi_n)-V(x_0)& = \la x_0 + \psi_n,\hat{B}_n (x_0+\psi_n)\ra - \la x_0,{B_0} x_0 \ra \\
 &\leq \la x_0 + \psi_n,\hat{B}_n (x_0+\psi_n)\ra - \la x_0,\hat{B}_n x_0 \ra.
\end{align*}
Again by the Mean Value Theorem,
\begin{equation*}
\begin{split}
& |\la x_0 + \psi_n,\hat{B}_n (x_0+\psi_n)\ra - \la x_0,\hat{B}_n x_0 \ra - 2 \la\psi_n,\hat{B}_n x_0\ra|\\
& \qquad\qquad\qquad\qquad\qquad\qquad\leq \|\psi_n\|_X \int_0^1 \Phi(\hat{B}_n,x_0+\psi_n)\,dt \leq \frac\varepsilon2 \|\psi_n\|_X
\end{split}
\end{equation*}
for all $n$ large enough. Thus,
\begin{equation*}
 V(x_0+\psi_n)-V(x_0) \leq 2 \la\psi_n,\hat{B}_n x_0\ra + \frac\varepsilon2 \|\psi_n\|_X.
\end{equation*}
Using one more time the compactness of $\B$, there exist a subsequence $n(k)$
such that $\hat{B}_{n(k)}\wot\hat{B}$ for some $\hat{B}\in\B$. 
Moreover, 
by continuity of $V$ and because of \eqref{wstard}, $\hat{B}\in \Scal(x_0)$. 
Hence,
\begin{align*}
 V(x_0+\psi_{n(k)})-V(x_0) &\leq 2 \la\psi_{n(k)},\hat{B} x_0\ra + 2 \la\psi_{n(k)},(\hat{B}_{n(k)}-\hat{B}) x_0\ra + \frac\varepsilon2 \|\psi_{n(k)}\|_X\\
&\leq \max_{B\in\Scal(x_0)} 2 \la \psi_{n(k)},Bx_0 \ra + \frac{3}{4}\varepsilon\|\psi_{n(k)}\|
\end{align*}
where we used that $2\la\psi_{n(k)},(\hat{B}_{n(k)}-\hat{B}) x_0\ra \leq \frac\varepsilon4\|\psi_{n(k)}\|_X$ for $k$ sufficiently large. This contradicts \eqref{eq:proofcontra3} and completes the proof.
\end{proof}

%

The following Corollary is now a direct consequence of Lemma~\ref{lem:Hilbert} and the choice
of the Lyapunov function $V(\cdot)$ in the proof of Theorem~\ref{thm:VnonComparable}.

\begin{corollary}\label{corl:Hilbert}
 Let $X$ be a separable Hilbert space 
and assume that there exist constants $K \geq 1$ and $\mu>0$ such that
\begin{equation}\label{corl:GUES}
 \|T_{\sigma(\cdot)}(t)\|_{\Lcal(X)} \leq Ke^{-\mu t}, \quad t \geq 0,~\sigma(\cdot)\text{-uniformly}.
\end{equation}
Then, there exists a 
function 
$V\: X \to [0,\infty)$ 
such that
$\sqrt{V(\cdot)}$ is a norm on $X$, $V(\cdot)$ is directionally Fr\'{e}chet differentiable,
\begin{equation}\label{corl:VBounds}
	c\|x\|_X^2 \leq V(x) \leq C \|x\|_X^2, \quad x\in X,
\end{equation}
for constants $c,C>0$ and
\begin{equation}\label{corl:LjVest}
	\Linf_j V(x) \leq -\|x\|_X^2, \quad j \in Q,~x\in X,
\end{equation}
with $\Linf_j V(x)$ defined as in \eqref{eq:LinfjVdef}.
\end{corollary}

\section{FINAL REMARKS AND OPEN PROBLEMS}\label{sec:conclusions}

We presented necessary and sufficient conditions for a (possibly infinite) family of semigroups to be 
globally uniformly exponentially stable for arbitrary switching signals $\sigma(\cdot)$, in terms of the 
existence of a common Lyapunov function. In particular, our results apply to switched  dynamical 
systems such as \eqref{eq:ssOpForm}, involving (possibly unbounded) operators on a Banach space $X$.

%


We have shown that the existence of a  norm, 
decaying uniformly along trajectories
and either 
bounded from above by a multiple of the Banach norm 
in presence of 
a uniform 
exponential growth bound for the switched system or 
comparable with the Banach norm
is equivalent to the switched system being globally uniformly exponentially stable. The latter  
generalizes a 
well-known
result for switched linear 
dynamical systems in 
$\RR^n$, $n\in\NN$ \cite{MolchanovPyatnitskiy1986}.
In the case in which $X$ is a separable Hilbert space the common Lyapunov function is shown to be Fr\'echet directionally differentiable. 

As an application, our results answer for example the question of existence of a common Lyapunov 
function for the switched linear hyperbolic system with reflecting boundaries considered in \cite{AminHanteBayen2008}.


Concerning the differences between Theorems~\ref{thm:VnonComparable} and~\ref{thm:VComparable}, we already noticed that a Lyapunov function $V(\cdot)$ satisfying condition ii) of the 
statement of Theorem~\ref{thm:VnonComparable} does not necessarily satisfy the stronger conditions appearing in the statement of Theorem~\ref{thm:VComparable}. Hence, Theorem~\ref{thm:VnonComparable} is better suited for proving 
the global uniform exponential stability of a switched system (although the uniform exponential growth boundedness should also be proved), while Theorem~\ref{thm:VComparable} provides more information on a switched system that is known to be globally uniformly exponentially stable, by tightening the properties satisfied by $V(\cdot)$. 

In the same spirit one could characterize 
the global uniform exponential stability of a switched system by 
(a priori) loosening  the hypotheses on $V(\cdot)$, 
replacing inequalities \eqref{thm1:LjVest} and \eqref{thm2:LjVest} 
by 
$$\Linf_j V(x)\leq -\kappa \|x\|_X^2, \quad\text{for all}~j \in Q,~\text{for some}~\kappa>0.$$ 
One could also remove the hypothesis that 
$V$ is the square of a norm. 
Conversely,  one could (a priori) tighten the same hypotheses replacing $\Linf_j V(x)$ by
$$\Lsup_j V(x)=\limsup_{t \downarrow 0} \frac{V(T_j(t)x) - V(x)}{t}.$$ 

%
As an open problem it remains 
to understand if smoothing can improve the 
regularity properties of the Lyapunov function $V(\cdot)$. For the 
finite-dimensional case, it is for instance known that $V(\cdot)$ can be taken polyhedral or polynomial 
\cite{Blanchini1994,BlanchiniMiani1999,DayawansaMartin1999}.
It would be interesting to recover results in the same direction for 
infinite dimensional switched systems.

\bibliographystyle{plain}
\bibliography{convTHMs}

\end{document}